\documentclass{article}

\usepackage[english]{babel}

\usepackage[letterpaper,top=2cm,bottom=2cm,left=3cm,right=3cm,marginparwidth=1.75cm]{geometry}

\usepackage{amsmath}
\usepackage{graphicx}
\usepackage[colorlinks=true, allcolors=blue]{hyperref}

\newcommand{\bl}{\color{black}}

\newtheorem{proposition}{Proposition}

\newtheorem{lemma}{Lemma}

\newtheorem{theorem}{Theorem}

\newtheorem{remark}{Remark}

\title{
New Lyapunov functions for
systems with source terms
\\
\emph{Dedicated to Jan Sokolowski}
}

\author{Martin Gugat
\thanks{Friedrich-Alexander-Universit\"at Erlangen-N\"urnberg (FAU),
Department Mathematik,
Lehrstuhl f\"ur Dynamics, Control, Machine Learning and Numerics (Alexander von Humboldt-Professur), 
Cauerstr. 11,
91058 Erlangen, Germany, 
({\tt martin.gugat@fau.de})}
}

\begin{document}
\maketitle

\begin{abstract}
Lyapunov functions with exponential weights 
have been used successfully 
as a powerful tool
for the stability analysis of hyperbolic systems of balance laws.
In this paper we extend the class of weight functions
to a family  of hyperbolic functions and
study the advantages in the analysis of
$2\times 2$ systems of balance laws.
We present cases
connected with the study of the limit of stabilizability 
where the new weights 
provide Lyapunov functions that show exponential
stability for a larger set of problem parameters than
classical exponential weights.

Moreover, we show that sufficiently large  time-delays
influence the limit of stabilizability in the sense that 
the parameter set where the system can be stabilized 
becomes substantially smaller.

We also demonstrate that  the hyperbolic weights are useful in
the analysis of the boundary feedback stability of systems of balance laws
that are governed by quasilinear hyperbolic  partial differential equations.
\end{abstract}

\section{Introduction}

In 
\cite{zbMATH01367198},
\cite{zbMATH06752608}
and related work 
exponential weights in a quadratic function
have been used 
to obtain a strict Lyapunov function for the stabilization 
of the Euler equation of incompressible fluids.
This valuable tool has been 
the key to achieve numerous
stabilization results for 
various systems,
see for example the survey paper \cite{MR4348713}.
In 
\cite{bastin2011boundary}
basic quadratic control Lyapunov function for linearized systems 
are investigated.
In  \cite{hayat2021global} 
also  nonlocal   source terms are studied for semilinear systems.

In many engineering applications 
that involve 
systems which 
can be modelled by hyperbolic systems  of partial differential
equations 
the source terms in these equations play an essential role.
In order to adapt the Lyapunov function candidates to this
situation, in this paper
we 
extend the exponential weights by
introducing a family of hyperbolic weight functions.
In the analysis, these Lyapunov functions
yield additional terms that help to obtain bounds for
the size of the admissible source terms, for which the system can
be stabilized from the boundary.

To illustrate the results that can be achieved with the hyperbolic weight functions,
we look at an example from
\cite{MR3561145}
that illustrates the limits of  boundary stabilizability,
see also
\cite{zbMATH07117070}.
We show that the hyperbolic  weights allow to extend the 
set where the Lyapunov function can be used to prove  the exponential stability.
The related  numerical  aspects of the boundary feedback stabilization for semilinear hyperbolic systems  have been studied in \cite{gerster2023numerical}.
The limits of stabilization of a networked linear hyperbolic system
with a circle have been studied in \cite{gugat2023limits}.
 
This paper has the following structure.
In Section \ref{examplebyBastinandCoron}
we present the  example by Bastin and Coron for the limits of stabilizability.
Section \ref{defiandprop} contains the definition and properties of the hyperbolic weight functions.
In Section \ref{sufficientcondifor}
we show that Lyapunov functions with the hyperbolic weights
yield sharper sufficient condition for stabilizability than the exponential  weights.
In Section \ref{linearweights} we use Lyapunov functions with linear weights to show
that for feedback gains that are too large, the system becomes unstable.
In Section \ref{delay} we study the influence of time delay in the boundary feedback on the stabilizability:
We show that 
for a time-delay that is sufficiently large, the system
becomes unstable even if it would be stable without time delay.
In Section  \ref{quasilinear} we show that Lyapunov functions with 
hyperbolic weights are also useful for the study of the stability of quasilinear systems
In Section \ref{quasilinearin} we present sufficient conditions
for the instability of systems of balance laws for sufiiciently long space intervals.




\section{The example by Bastin and Coron}
\label{examplebyBastinandCoron}

\noindent
In \cite{MR3561145},
Bastin and Coron consider the following  system in diagonal form:
\begin{align}
\label{bastincoron1}
(\delta_+)_t + (\delta_+)_x +  {\cal M}\, \delta_- & =  0,
\\
\label{bastincoron2}
(\delta_-)_t - (\delta_-)_x +  {\cal M} \,  \delta_+ & = 0.
\end{align}
Here  $ {\cal M} >0$ is  a real parameter,
 $x$ is in the interval $(0,\, L)$ and
$t\geq 0$.
With the boundary conditions 
\begin{equation}
\label{riemannfeedback0}
\delta_+(t,\, 0) = k \, \delta_-(t,\, 0),
\end{equation}
\begin{equation}
\label{rb15062018a}
\delta_-(t,\, L)= \delta_+(t,\, L)
\end{equation}
and initial states
$\delta_+(0, \cdot), \delta_-(0, \, \cdot)
\in H^1(0, \, L)$
the system is completed.
The following proposition is shown:
\begin{proposition}
	\label{proposition1}
	If 
 \[   {\cal M} \, L \geq \pi,\]
 there is no real value of $k$ such that the closed loop system
	(\ref{bastincoron1}), (\ref{bastincoron2}),
	(\ref{riemannfeedback0}),
	(\ref{rb15062018a})
	is exponentially stable.
\end{proposition}

%

%
%
%
%
%
%
%


\subsection{The example by Bastin and Coron: 
\\
A sufficient condition for stabilizability
with exponential weights}

The following Proposition from \cite{zbMATH07117070} 
is proved using a Lyapunov function with exponential weights.
It states that if $L>0$ is sufficiently small,
the closed loop system is exponentially stable  if $|k|<1$.

\begin{proposition}
	\label{proposition1a}
	If with $\lambda>0$ we have
	$|k|\leq {\rm e}^{-\lambda\, L}$
	and
	$ {\cal M}     < \tfrac{\lambda}{  1 + {\rm e}^{2 \, \lambda \, L} } $,
	the
	closed loop system
	(\ref{bastincoron1}), (\ref{bastincoron2}),
	(\ref{riemannfeedback0}),
	(\ref{rb15062018a})
	is exponentially stable
	for all initial states
	$(\delta_+(0, \cdot), \delta_-(0,  \cdot))\in (H^1(0,\, L))^2$.
\end{proposition}
The proof is presented in \cite{zbMATH07117070} 
using  the Lyapunov function
	\[L(t) = \frac{1}{2}\int_0^L \exp(\lambda(L-x)) \, \delta_+^2(t,\, x) + \exp(\lambda(x-L)) \, \delta_-^2(t,\, x)\, dx.\]
 Note that the notation in \cite{zbMATH07117070}  is different,
 namely $(U, \, V)$ instead of $(\delta_+, \delta_-)$.

 Proposition  \ref{proposition1a} yields stability only if the following inequality holds:
\begin{equation}
\label{alteschranke}
 {\cal M} L< \sup_{z>0} \tfrac{z}{  1 + {\rm e}^{2 \, z} }  =  \frac{W(\exp(-1))}{2}  =  0.139...
 \end{equation}
 where $W$ is the Lambert W-function.



\section{Definition and properties of the hyperbolic weight functions}
\label{defiandprop}

In this section we define the hyperbolic weight functions that generalize the exponential
weights that have been used for example in
\cite{zbMATH06752608}.
The exponential weights 
$\exp(\mp \psi \, x)$ 
occur naturally since their derivatives 
can again be expressed in terms of these weights,
that is they satisfy a linear differential equation. 
Moreover, they can be used as weight functions  since they only attain positive values.

Since $\exp( \pm \psi x) = \cosh( \psi x) \pm \sinh( \psi x)$,
a natural perturbation is the weight function
\[
h_\pm(x) = \sqrt{\upsilon} \cosh( \psi \, x) \mp \sinh( \psi x)
\]
where $\upsilon >0$ is chosen sufficiently large such
that only positive values are attained.

We have the representation
\[h_\pm(x) =  \cosh(\psi \, x) \left[\sqrt{\upsilon} \mp \tanh(\psi \, x) \right]
.
\]
Thus  if $\sqrt{\upsilon}  > |\tanh(\psi \, L)|$, we have
$h_\pm(x) >0$ for all  $x\in [-L, \, L]$
and thus the functions $h_\pm(x)$ can be used as weight 
functions  in Lyapunov functions for subintervals of $[-L, L]$.
The derivatives of the  weight functions
can be represented as a linear
combination of  the  weight functions.

\begin{figure}
    \centering
        \includegraphics[width=\textwidth]{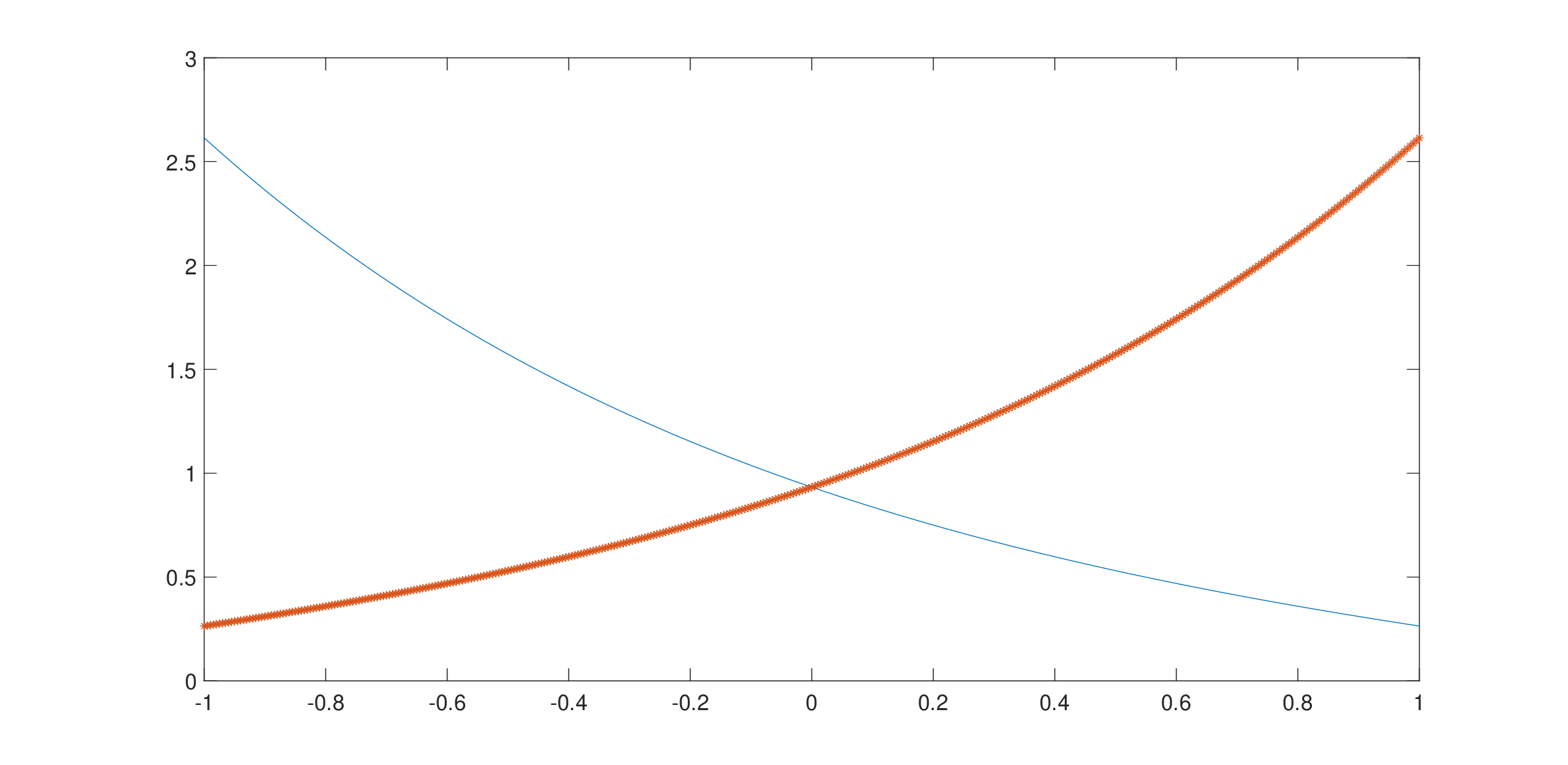}
        \caption{The weight functions $h_+$ and $h_-$ (with '$\ast$') 
        for $\psi = L = 1$ and $\nu =  \tfrac{3}{2} \,   \tanh^2(\psi \, L) $.}
        \label{fig:Ex1}
\end{figure}

In the following Lemma we summarize properties of 
$h_\pm(x)$, and state how
the derivatives can be expressed in terms of the weight functions.
Note that the exponential weights
\[
h_+(x) = \exp(- \psi \, x), \;\;
h_-(x) = \exp( \psi \, x)
\]
occur as the special case $\upsilon =1$.
For $\psi \rightarrow 0+$
both $h_+$ and $h_-$ converge to the constant function 
$\sqrt{\upsilon}$.
\\ 
Figure \ref{fig:Ex1} shows the graphs of
$h_+$ and $h_-$ for $\psi = L = 1$ and $\upsilon =  \tfrac{3}{2} \,   \tanh^2(\psi \, L) $.
We have $\sqrt{ \upsilon } = 0.9328... $

\begin{figure}
    \centering
        \includegraphics[width=\textwidth]{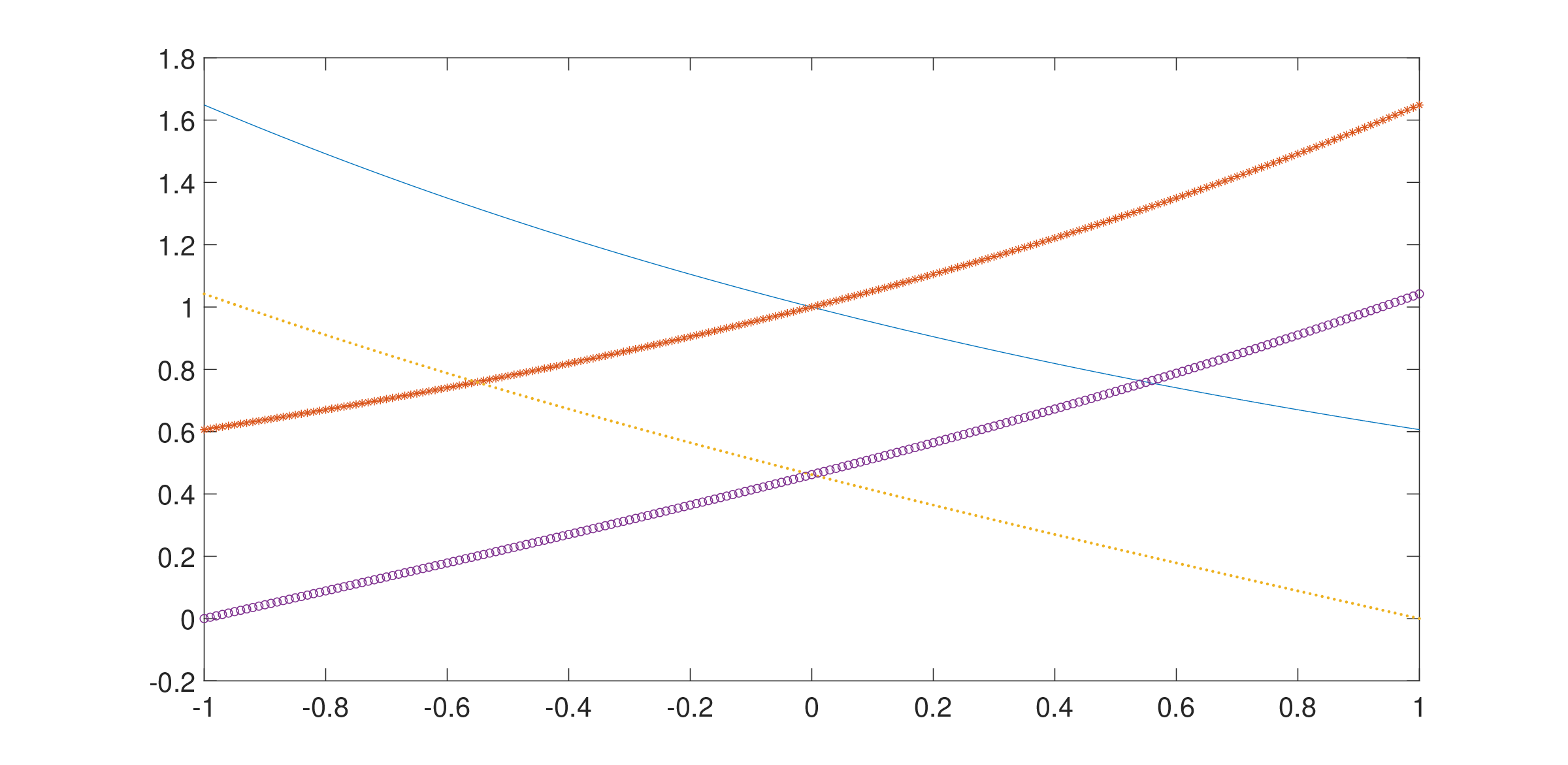}
        \caption{The weight functions $h_+$ and $h_-$ 
        for $L = 1$, $\psi = \tfrac{1}{2}  $, $\upsilon = 1$ and $\upsilon =  \tanh^2(\psi \, L) $.
      The case of the exponential weights is  $ \upsilon= 1$  with $h_+$: full line, $h_-: '*'$.
      In this case $h_+(0)= h_-(0)=1$. 
      The extremal case $\upsilon =   \tanh^2(\psi \, L) $ is shown 
  with $h_+$: 'o', $h_-$: $\cdot$. 
  In this case $h_+(1)= h_-(-1)=0$.
        }
        \label{fig:Ex2}
\end{figure}

Figure \ref{fig:Ex2} shows the graphs of
$h_+$ and $h_-$ for $L = 1$, $\psi = \frac{1}{2}$ and
$ \upsilon \in\{ 1,  \,  \tanh^2(\psi \, L)\}$.
The case $ \upsilon = 1$ is the case of the exponential weights and shown with 
$h_+$: full line, $h_-: '*'$.
The extremal case 
$ \upsilon  =   \tanh^2(\psi \, L) $
where the weight functions loose positivity since 
$h_+(L) = 0$
and
$h_-(-L)=0$
is shown with  $h_+$: 'o', $h_-$: $\cdot$.

\begin{lemma}
\label{lemma1}
Let $\psi>0$ and $\upsilon > \tanh^2(\psi \, L) $
be given.
For $x\in [-L, \, L]$, define the
functions 
\[
h_\pm(x) =
  \sqrt{\upsilon  }
 \cosh(\psi\, x) \mp \sinh(\psi\, x) .
\]

\textbf{a)} The functions $h_+$ and $h_-$ only attain values in $(0, \infty)$ for $x\in [- L, \, L]$.

Since $h''_\pm(x) = \psi^2 \,  h_\pm(x)$, this implies that
the functions $h_\pm(x)$ are strictly convex on $[-L, \, L]$.

\textbf{b)}
We have $h_+(0) = h_-(0) =  \sqrt{\upsilon}$.
We have
$h_+'(0) = - \psi$  and $h_-'(0) = \psi$.

\textbf{c)} 
We have
\[
\frac{h_+(x)}{h_-(x)}
=
\frac{2 \, \sqrt{\upsilon  }
}{ 
\sqrt{\upsilon  }
+ \tanh( \psi \, x)   } - 1.
\]

Hence  $\frac{d}{dx} \left( \frac{h_+(x)}{h_-(x)} \right) < 0 $. 
Thus   $\frac{h_+(x)}{h_-(x)}$ is decreasing
and
 $\frac{h_-(x)}{h_+(x)}$ is increasing. 
 
For $x\in (0, L)$ we have 
$0 <\frac{h_+(x)}{h_-(x)}< 
\frac{h_+(0)}{h_-(0)} = 1$
and
$
\frac{h_-(x)}{h_+(x)} 
< \frac{h_-(L)}{h_+(L)} 
=
\frac{ \upsilon^{\frac{1}{2}}   +   \tanh (\psi \, L)  }{ 
\upsilon^{\frac{1}{2}} - \tanh (\psi \, L)}  
$.
We have
\begin{equation}
    \label{quationt23022024}
\frac{h_+(\tfrac{L}{2})}{h_-(\tfrac{L}{2})} = \frac{ \sqrt{\upsilon} - \tanh(\psi \,\tfrac{L}{2} ) }{ \sqrt{\upsilon} + \tanh(\psi \, \tfrac{L}{2}) }
=
\frac{ \sqrt{\upsilon} \left(1 + \tfrac{1}{\cosh(\psi\, L)} \right) - \tanh(\psi \, L)}{ \sqrt{\upsilon} \left(1 + \tfrac{1}{\cosh(\psi\, L)} \right) +  \tanh(\psi \, L)}.
\end{equation}

\textbf{d)} 
For $x\in [-L, L]$ we have
\[h_-(x) = h_+(-x).\]

\textbf{e)}
We have
\[h_+'(x) =
- \psi \left[  \frac{   \upsilon  + 1 }{  2 \sqrt{\upsilon}   }      \right] h_+(x)
+ \psi \left[ \frac{ \upsilon  - 1  }{   2 \sqrt{\upsilon}     }  \right] h_-(x),
\]

\[h_-'(x) = - \psi \left[ \frac{ \upsilon  - 1  }{   2 \sqrt{\upsilon}     }  \right] 
 h_+(x)
 +
 \psi \left[  \frac{   \upsilon  + 1 }{  2 \sqrt{\upsilon}   }      \right]
 h_-(x).
\]

 For $\upsilon \in  (\tanh^2(\psi \, L) , 1]$,
 $h_+$ is decreasing and $h_-$ is increasing on $[-L, L]$.

\end{lemma}


\textbf{Proof:}

\textbf{a)}   
Since
$\sqrt{\upsilon} > \tanh(\psi \, L) $ for all $x\in [-L, \, L]$ we have
\[
h_\pm(x) =    \cosh(\psi\, x)  \left[ 
\sqrt{\upsilon  }
\mp \tanh(\psi\, x) \right] >0.
\]

\textbf{e)}
We have
\[h_\pm'(x) = 
\psi \, 
\sqrt{\upsilon  }
\sinh( \psi\, x) \mp  \psi \cosh(\psi\, x) .
\]
Since 
\[  \cosh(\psi\, x) =   \frac{1}{2\, \sqrt{\upsilon}  } \left[ h_+(x) + h_-(x)  \right],\;\;
\sinh(\psi\, x) = \frac{1}{2} \left[  - h_+(x) + h_-(x) \right]
\]
this yields
\[h_\pm'(x) = 
\psi \,
\sqrt{\upsilon  }
\frac{1}{2}
\left[    - h_+(x) + h_-(x)   \right]     \mp  \psi   \frac{1}{2\, \sqrt{\upsilon}  } \left[ h_+(x) + h_-(x)  \right].  
\]


Thus we have
\[
h_\pm'(x) =
- \psi \left[  \frac{   \upsilon  \pm 1 }{  2 \sqrt{\upsilon}   }      \right] h_+(x)
+ \psi \left[ \frac{ \upsilon  \mp  1  }{   2 \sqrt{\upsilon}     }  \right] h_-(x).
\]



Since
\[h_\pm'(x) = 
\mp \psi \, 
\sqrt{\upsilon  }
\cosh(\psi\, x) 
\left[ 
\frac{1}{ 
\sqrt{\upsilon  }
}
\pm 
 \tanh( \psi\, x) \right]
\]
and $ \frac{1}{ 
\sqrt{\upsilon  }
} \geq 1 $
we have
$h_+'(x) < 0$ for all $x\in [-L, \, L]$.
Moreover
we have
$h_-'(x) > 0$ for all $x\in [-L, \, L]$.

\textbf{b)} We have $ h_\pm(0) =    \cosh(\psi\, 0)  \left[
\sqrt{\upsilon  }
\mp \tanh(\psi\, 0)\ \right]   =
\sqrt{\upsilon  }
.  $
  We have
  $h_\pm'(0) =\mp \psi  $.

\textbf{c)} 
We have
\[
\frac{h_+(x)}{h_-(x)}
=
\frac{  \upsilon^{\frac{1}{2}}  - \tanh(\psi\, x)  }{  \upsilon^{\frac{1}{2}} + \tanh(\psi\, x)  }
=
\frac{2 \, \upsilon^{\tfrac{1}{2}}}{  \upsilon^{\tfrac{1}{2}} + \tanh( \psi \, x)   } - 1.
\]

\textbf{d)} 
For $x\in [-L, L]$ we have
\[ h_+(-x)  =
  \cosh(-\psi\, x)  \left[ \upsilon^{\frac{1}{2}}  - \tanh(-\psi\, x) \right] 
  =
  \cosh(\psi\, x)  \left[ \upsilon^{\frac{1}{2}}  + \tanh(\psi\, x) \right]  = h_-(x).
\]


\section{The example by Bastin and Coron: A sufficient condition for stabilizability
with hyperbolic weights}
\label{sufficientcondifor}

Let $\psi>0$ be given. 
Define the hyperbolic  weights $h_+(x)>0$, $h_-(x)>0$  ($x\in [-L , \, L]$)
such that for 
$\upsilon > \tanh^2(\psi \, L) $
 we have 
\[
h_\pm(x) =   
\upsilon^{\frac{1}{2}}   \cosh(\psi\, x) \mp \sinh(\psi\, x)
.
\]
For $\upsilon =1$ we obtain the exponential weights
$
h_\pm(x)   =      \exp( \mp \psi x)
$
.

%

Define the Lyapunov candidate  function
\begin{equation}
    \label{lyapdefi}
{\cal E}(t) := 
{\frac{1}{2}} 
\int_0^{L}
h_+(x - L) \, |\delta_+(t,\, x)|^2 + 
h_-(x - L) \, |\delta_-(t,\, x)|^2 \, dx.
\end{equation}
	For the time-derivative we obtain using (\ref{bastincoron1}) and  (\ref{bastincoron2}) 
	\begin{align*}
		{\cal E}'(t) & =  \int_0^L h_+(x- L) \, \delta_+\, (\delta_+)_t + 
   h_-(x- L)\, \delta_-\, (\delta_-)_t\, dx
		\\
		& = 
		\int_0^L h_+(x- L) \, \delta_+\, (-  (\delta_+)_x -   {\cal M}   \,  \delta_-) 
		\quad +   h_-(x- L ) \, \delta_-\, (  (\delta_-)_x -    {\cal M} \, \delta_+)\, dx
		\\
		& = 
		\int_0^L    -   h_+(x- L)
                               \left(\tfrac{1}{2} (\delta_+)^2 \right)_x 
		 +
		h_-(x- L)\, \left(\tfrac{1}{2} (\delta_-)^2 \right)_x
		-   {\cal M} \left[h_+(x-L) + h_-(x-L)\right]   \, \delta_+\, \delta_- \, dx.
	\end{align*}

	Integration by parts yields
	\begin{align*}
	        {\cal E}'(t)     = 
		&  \tfrac{1}{2} \int_0^L  
     h_+'(x-L) \, \delta_+^2 -  h_-'(x-L) \, (\delta_-)^2 
  	-  2 {\cal M} \left[  h_+(x-L) + h_-(x-L)  \right]   \, \delta_+\, \delta_- \, dx 
     \\
		&+
		\left[\tfrac{1}{2}  h_-(x-L) \, (\delta_-(t,x))^2  
		-   \tfrac{1}{2} h_+(x-L) \, (\delta_+(t,x))^2\right] \Big|_{x=0}^{x=L}
  \\
  & = - \tfrac{\psi}{2}  \int_0^L \left[  \frac{   \upsilon  + 1 }{  2 \sqrt{\upsilon}   }      \right] h_+(x-L)\, \delta_+^2
  + \left[ \frac{ 1 - \upsilon   }{   2 \sqrt{\upsilon}     }  \right] h_-(x-L) \, \delta_+^2
   \\ &  + \left[ \frac{ 1 - \upsilon  }{   2 \sqrt{\upsilon}     }  \right] 
 h_+(x - L)  \, \delta_-^2 
  + \left[  \frac{   \upsilon  + 1 }{  2 \sqrt{\upsilon}   }      \right]
 h_-(x-L)   \, \delta_-^2 \, dx 
 \\ &  - \int_0^L  {\cal M} \left[  h_+(x-L) + h_-(x-L)  \right]   \, \delta_+\, \delta_- \, dx
  \\
		&+
		\left[\tfrac{1}{2}  h_-(x-L) \, (\delta_-(t,x))^2  
		-   \tfrac{1}{2} h_+(x-L) \, (\delta_+(t,x))^2\right] \Big|_{x=0}^{x=L}.
	\end{align*}
 If the boundary term in the last line is less than or equal to zero, 
  for $\upsilon \in ( \tanh^2(\psi \, L)   , 1]$
  using $|z_1\, z_2|  \leq \frac{|z_1|^2}{2} + \frac{|z_2|^2}{2}$   we obtain the inequality 
\begin{align*}
    {\cal E}'(t)    \leq 
  & - \tfrac{\psi}{2}  \int_0^L \left[  \frac{   \upsilon  + 1 }{  2 \sqrt{\upsilon}   }  -    \frac{{\cal M}}{\psi}
  \right] \, h_+(x - L)\, \delta_+^2
  + \left[ \frac{ 1 - \upsilon   }{   2 \sqrt{\upsilon}     } 
   -      \frac{{\cal M}}{\psi}
  \right] \, h_-(x - L) \, \delta_+^2
   \\ &  + \left[ \frac{ 1 - \upsilon  }{   2 \sqrt{\upsilon}     }  -     \frac{{\cal M}}{\psi}  \right] \,
 h_+(x- L)  \, \delta_-^2 
  + \left[  \frac{   \upsilon  + 1 }{  2 \sqrt{\upsilon}   }    -    \frac{{\cal M}}{\psi}  \right] \, 
 h_-(x- L)   \, \delta_-^2 \, dx. 
	\end{align*}
If
\begin{equation}
\label{10072023}
{\cal M} 
\leq
\frac{\psi}{2}  \, 
\frac{ 1 - \upsilon  }{    \sqrt{\upsilon}     }
=
\frac{\psi}{2}  \, 
\left(
\frac{1}{    \sqrt{\upsilon}     }
-
{    \sqrt{\upsilon}     }
\right)
\end{equation}
due to   (\ref{lyapdefi})    this yields
\begin{align*}
    {\cal E}'(t)    \leq  - 
\psi
\left[  \frac{   \upsilon  + 1 }{  2 \sqrt{\upsilon}   }   -    \frac{ 1 - \upsilon  }{   2 \sqrt{\upsilon}     } \right]    {\cal E}(t) 
=
 - \psi
\,
\sqrt{\upsilon} 
\,
{\cal E}(t).
	\end{align*}
	By Gronwall's inequality this implies
	$ {\cal E}(t) \leq \exp(-  
 \psi 
\,
\sqrt{\upsilon}  \, t) \, 
{\cal E}(0),
$
	so it remains to check the negativity assumption on
 the boundary term.
 Note that (\ref{10072023}) is equivalent to
\begin{equation}
\label{10072023xyz}
{\cal M} L
\leq
\frac{\psi \, L}{2} 
\frac{ 1 - \upsilon  }{    \sqrt{\upsilon}     }
=
\frac{\psi L}{2} 
\left(
\frac{1}{    \sqrt{\upsilon}    }
-
{    \sqrt{\upsilon}     }
\right).
\end{equation}
For the right--hand side
of (\ref{10072023xyz}) 
we have the upper  bound
\[
\frac{\psi L}{2} \, 
\left(
\frac{1}{    \sqrt{\upsilon}    }
-
{    \sqrt{\upsilon}     }
\right)
\leq
\frac{ 1 }{2} \, 
\sup_{z>0}
z \left(
\frac{1}{  \tanh(z)   }
-
{   \tanh(z)    }
\right)
=
\frac{ 1 }{2} \, 
\sup_{z>0}
\frac{z}{  \cosh(z) \, \sinh(z) }
\]
\[
=
2 \,
\sup_{z>0}
\frac{ z}{ \exp(2\, z) - \exp(- 2\, z) }
=
2 \,
\sup_{z>0}
\frac{ z}{ 4\, z +  \frac{16}{6} z^3 + .....  }
= \frac{1}{2}.
\]
So here we obtain an upper bound for
${\cal M} L$ that is
 closer to $\pi $ than the
value 
from  (\ref{alteschranke}) 
that is obtained with the exponential weights.

Since $h_+(0) = h_-(0)$,
for the boundary term at  $x=L$ we get
	\[
	\tfrac{1}{2} (\delta_-)^2(t,\, L)  -  \tfrac{1}{2} (\delta_+)^2(t,\, L)  = 0.
	\]

For the boundary term at  $x=0$ we get
	\[
	 h_+(-L) \, k^2   -  h_-(-L) 
	\leq 0
 \]
 if and only if
 \[
  k^2  \leq \frac{h_-(-L)  }{ h_+(-L) }
  = 
   \frac{h_+(L)  }{ h_-(L) }.
 \]
 This is the case it
$k^2 \leq  \frac{ \sqrt{\upsilon} - \tanh(\psi \, L)  }{ \sqrt{\upsilon} + \tanh(\psi L) }$.
 %
%
 To be precise, in this case it suffices that
    \[
    k^2 <  B:= \sup_{\psi >0} \sup_{ \upsilon \in ( \tanh^2(\psi \, L)   , 1)    }  \frac{ \sqrt{\upsilon} - \tanh(\psi \, L)  }{ \sqrt{\upsilon} + \tanh(\psi L) }
    .
    \]
        With the choice $\sqrt{\upsilon} = \frac{1 + \tanh(\psi \, L)}{2}$ we have 
       $ \frac{ \sqrt{\upsilon} - \tanh(\psi \, L)  }{ \sqrt{\upsilon} + \tanh(\psi L) }
        = \frac{    1 - \tanh(\psi \, L)  }{ 1  + 3\tanh(\psi L) }$.
\\
This yields $B  \geq \lim\limits_{\psi \rightarrow 0 +}  \frac{    1 - \tanh(\psi \, L)  }{ 1  + 3\tanh(\psi L) } = 1$.

Thus using the Lyapunov function with hyperbolic weights 
have shown the following result:
\begin{proposition}
	\label{proposition1b}
	If   
 \[{\cal M} L<  \frac{1}{2}  \]
        and 
	$|k|$ is sufficiently small
 (in the sense that  $|k| < 1$) 
	the	closed loop system
	(\ref{bastincoron1}), (\ref{bastincoron2}),
	(\ref{riemannfeedback0}),
	(\ref{rb15062018a})
	is exponentially stable
	for all initial states
	$(\delta_+(0, \cdot), \delta_-(0, \cdot))\in (H^1(0,\, L))^2$
        with $\delta_+(0, \, L) = \delta_-(0, L) $
        and $\delta_+(0,0) = k \, \delta_-(0,0)$.

\end{proposition}

Note that the stability result
in Proposition \ref{proposition1b}
holds for a larger set of parameters than 
the bound (\ref{alteschranke}) that is implied by 
Proposition \ref{proposition1a}.


\section{The example by Bastin and Coron continued:
\\
A sufficient condition for instability
with affine linear  weights}
\label{linearweights}
 
In this section introduce  a Lyapunov function with affine linear weights to show 
 that the system is unstable if $|k|$ is too large.
This illustrates further  the flexibility of the analysis 
that is based upon Lyapunov functions.
Similarly as for the hyperbolic weights,
the derivatives of the affine linear weights,
can be represented as a linear
combination of affine linear weights.

%
Define the affine linear 
weights $h_\pm(x)>0$
by
\[
h_\pm(x) =   
1 \pm   
2
\,
{\cal M}\,  x
.
\]


Note that in the case without source term, that is with
 ${\cal M}=0$ the definition yields  $h_\pm(x) = 1$,
 that is constant weights.
 If
$  {\cal M}\, L <  \frac{1}{2}$
we have $h_\pm(x) >0$ for all $x\in [-L , \, L]$.
We have
\[h_+'(x) = 
{{\cal M}}
(h_+(x) + h_-(x)),\;h_-'(x) =  -  
{{\cal M}}
(h_+(x) + h_-(x)).
\]

Consider again  the Lyapunov candidate  function
\[
{\cal E}(t) :=  \frac{1}{2} \int_0^{L}
h_+(x - L) \, |\delta_+(t,\, x)|^2 + 
h_-(x - L) \, |\delta_-(t,\, x)|^2 \, dx.
\]
	For the time-derivative we obtain as above using integration by parts 
%
	\begin{align*}
	        {\cal E}'(t)     = 
		&  \tfrac{1}{2} \int_0^L  
     h_+'(x-L) \, \delta_+^2 -  h_-'(x-L) \, \delta_-^2 
  	-  2 {\cal M} \left[  h_+(x-L) + h_-(x-L)  \right]   \, \delta_+\, \delta_- \, dx
     \\
		&+
		\left[\tfrac{1}{2}  h_-(x-L) \, (\delta_-(t,x))^2  
		-   \tfrac{1}{2} h_+(x-L) \, (\delta_+(t,x))^2\right] \Big|_{x=0}^{x=L}
  \\
  & =  \int_0^L  \frac{{\cal M}}{2} \,  h_+(x-L)\, \delta_+^2
  +    \frac{{\cal M}}{2} \, h_-(x-L) \, \delta_+^2
 +   \frac{{\cal M}}{2} h_+(x - L)  \, \delta_-^2 
  + 
  \frac{{\cal M}}{2}   h_-(x-L)   \, \delta_-^2 \, dx 
 \\ &  - \int_0^L  {\cal M} \left[  h_+(x-L) + h_-(x-L)  \right]   \, \delta_+\, \delta_- \, dx
  \\
		&+
		\left[\tfrac{1}{2}  h_-(x-L) \, (\delta_-(t,x))^2  
		-   \tfrac{1}{2} h_+(x-L) \, (\delta_+(t,x))^2\right] \Big|_{x=0}^{x=L}.
	\end{align*}
 We obtain the equation 
\begin{align*}
    {\cal E}'(t)     
  & =  \int_0^L   \frac{{\cal M}}{2}  \,  h_+(x-L)\, ( \delta_+ - \delta_-)^2
  +    \frac{{\cal M}}{2}  \, h_-(x-L) (\delta_+ - \delta_-)^2\, dx
  \\
		&+
		\left[\tfrac{1}{2}  h_-(x-L) \, (\delta_-(t,x))^2  
		-   \tfrac{1}{2} h_+(x-L) \, (\delta_+(t,x))^2\right] \Big|_{x=0}^{x=L}.
	\end{align*}
 
This yields
\begin{align*}
    {\cal E}'(t)     \geq 
		&
		\left[\tfrac{1}{2}  h_-(x-L) \, (\delta_-(t,x))^2  
		-   \tfrac{1}{2} h_+(x-L) \, (\delta_+(t,x))^2\right] \Big|_{x=0}^{x=L}.
	\end{align*}

Since $h_+(0) = h_-(0)$,
for the boundary term at  $x=L$ we get
	\[
	\tfrac{1}{2} (\delta_-)^2(t,\, L)  -  \tfrac{1}{2} (\delta_+)^2(t,\, L)  = 0.
	\]

For the boundary term at  $x=0$ we get
	\[
	 h_+(-L) \, k^2   -  h_-(-L) 
	\geq 0
 \]
 if and only if
 \[
  k^2  \geq \frac{h_-(-L)  }{ h_+(-L) }
  = 
   \frac{h_+(L)  }{ h_-(L) }.
 \]
 This is the case it
$k^2 \geq  \frac{ 1 + 2 {\cal M}\,  L  }{ 1 - 2 {\cal M}\,  L }$.
Then we have  ${\cal E}'(t)    \geq 0$.
 %
%
Thus using the Lyapunov function with affine linear  weights 
have shown the following result:
\begin{proposition}
	\label{proposition1bb}
	If   
 \[ {\cal M} L<  \frac{1}{2} \]
        and 
	
 \[
 k^2 \geq  \frac{ 1 + 2\, {\cal M}\,  L  }{ 1 - 2\,{\cal M}\,  L }
 \]
	the	closed loop system
	(\ref{bastincoron1}), (\ref{bastincoron2}),
	(\ref{riemannfeedback0}),
	(\ref{rb15062018a})
	is unstable
	for all initial states
	$(\delta_+(t,0), \delta_-(t,0))\in (H^1(0,\, L))^2$.
\end{proposition}

\begin{remark}
Note that in the case without source term
(that is ${\cal M}=0$) 
the proof also works and yields
instability for all
$|k|\geq 1$ (see Theorem 2.4 in  \cite{MR3561145}).

As pointed in \cite{MR3561145}),
the results from \cite{lichtner2008spectral} imply that
exponential stabilization can only be achieved 
if $|k|<1$.
So the bound provided in Proposition  \ref{proposition1bb}
is no novelty. The novelty 
is the construction of the Lypunov function
with affine linear weights that is used
for an easy proof of the statement.



\end{remark}



\section{The influence of time delay in the boundary feedback on the stabilizability}
\label{delay}

In this section
we discuss 
 the influence of time delay on the stabilizability of
the example by Bastin and Coron, that is 
the	closed loop system
	(\ref{bastincoron1}), (\ref{bastincoron2}),
	(\ref{riemannfeedback0}),
	(\ref{rb15062018a}).
In particular, we want to know whether a sufficiently
large time delay can lead to
non-stabilizability for a system that is stabilizable 
in the case without time-delay.
So we ask the question:
Can a sufficiently large time-delay lead to
a decrease of the critical length?
At this point, it is appropriate to mention
Datko's classical contributions to the study of 
time-delay, see \cite{datko1988not}
and \cite{datko1986example}
where it is shown that arbitrarily small time-delay
can destabilize a system that is otherwise stable.
A recent contribution on the topic 
for nonlinear systems is
\cite{haidar2021lyapunov}.
Our result is of a different type:
We show that if the time-delay is sufficiently large,
it 
can make the region where stabilization is possible substantially smaller.
So time-delay influences the limits of stabilizability.
This result is related to  \cite{gugat2011example},
where it is shown that  for certain time delays 
appropriately chosen sufficiently small feedback gains lead
to stability.

Let $\tau>0$ be a given time delay.
For $t\geq 2\, \tau$, we replace
the feedback law 
(\ref{riemannfeedback0})
by
\begin{equation}
\label{riemannfeedback0delay}
\delta_+(t,\, 0) = k \, \delta_-(t-\tau,\, 0).
\end{equation} 
To complete the system,
similarly as in \cite{MR2871937}
a compatible starting phase on
the time interval $[0, 2\,\tau)$
has to be added,
such that a well-defined regular system state is generated.
For $t\in [0, \,2\tau)$ we define
\[
\delta_+(t,\, 0) =  k \, \delta_-(t- \zeta(t),\, 0).
\]
We choose $\zeta(t)$ as a smooth function with
$\zeta(0)= 0$, $\zeta(2\,\tau) = \tau$, $\zeta'(2\tau)= 0$
and $\zeta(t) \leq t $ for all $t\in [0, 2\tau]$.
%
The following Proposition  provides an affirmative answer to the question posed above: 
\begin{proposition}
	\label{proposition1delay}
 Let  $\hat k > 0$ be given. 	If 
\begin{equation}
\label{assumption09082023}
{\cal M} L 
\in \left(\frac{3}{4} \pi, \, \pi \right) 
\end{equation}
and  $\tau>0$ is sufficiently large,
 there is no value of $k \in (-   \hat k, \,  \hat k ) $ such that the closed loop system
	(\ref{bastincoron1}), (\ref{bastincoron2}),
	(\ref{riemannfeedback0delay}),
	(\ref{rb15062018a})
	is exponentially stable.
\end{proposition}

\textbf{Proof:} 
To represent the state for $t \geq 2\tau$   we consider the separation ansatz
\[
\delta_+(t,\, x)  =  \exp(\sigma\, t) \, f(x),
\;\;
\delta_-(t,\, x)  =  \exp(\sigma\, t) \,  g(x)
\]
for the solution of (\ref{bastincoron1}), (\ref{bastincoron2}). 
For $\sigma\in (0,\, {\cal M})$ define
$
\omega = \sqrt{{\cal M}^2 - \sigma^2}>0.
$
The pdes (\ref{bastincoron1}) and  (\ref{bastincoron2}) imply
$f''(x) + ({\cal M}^2 - \sigma^2) f(x) = 0$,
and $g''(x)  + ({\cal M}^2 - \sigma^2) g(x) = 0$.

Hence  the solutions have the form
$f(x) = A \, \sin(\omega\, x) + B\, \cos(\omega\, x) $,
$g(x) = C \, \sin(\omega\, x) + D\, \cos(\omega\, x) $
with real numbers $A$, $B$, $C$, $D$. 
The pdes (\ref{bastincoron1}) and  (\ref{bastincoron2}) imply
\begin{equation}
\label{lgs}
\begin{pmatrix}
\sigma & - \omega & {\cal M} & 0
\\
 \omega &  \sigma & 0 &  {\cal M} 
 \\
  {\cal M} & 0 & \sigma &  \omega 
  \\
  0 &  {\cal M} & -\omega & \sigma
\end{pmatrix}
\,
\begin{pmatrix}
A
\\
B
\\
C
\\
D
\end{pmatrix}
=
\begin{pmatrix}
0
\\
0
\\
0
\\
0
\end{pmatrix}
.
\end{equation}
The feedback law 
(\ref{riemannfeedback0delay})
implies 
$B= f(0) = k \,  \exp( - \sigma \, \tau)  \, g(0) = 
  k \,  \exp( - \sigma \, \tau) \, D$.
We set $D = - \omega$.
Then $B =  - k \exp(-\sigma \, \tau) \,  \omega$ and
(\ref{lgs}) yields
\begin{equation}
\label{lgsstrich}
\begin{pmatrix}
\sigma &  {\cal M} 
\\
 \omega &  0 
\end{pmatrix}
\,
\begin{pmatrix}
A
\\
C
\end{pmatrix}
=
\begin{pmatrix}
\omega \, B
\\
-\,\sigma B + {\cal M} D
\end{pmatrix}
=
\begin{pmatrix}
- k \exp(-\sigma \, \tau) \,  \omega^2
\\
 k \exp(-\sigma \, \tau) \, \sigma\, \omega
 + {\cal M} \omega
\end{pmatrix}
.
\end{equation}
  This yields
  $ A=  k \exp(-\sigma \, \tau) \, \sigma
 + {\cal M}$
 and 
$ C = - k \exp(-\sigma \, \tau) \, {\cal M} - \sigma$.
Note that
\begin{equation}
\label{lgsstrich2}
\begin{pmatrix}
 {\cal M} & \sigma &
\\
0 & -\omega 
\end{pmatrix}
\,
\begin{pmatrix}
A
\\
C
\end{pmatrix}
=
\begin{pmatrix}
- \omega \, D
\\
- {\cal M} B  - \sigma D
\end{pmatrix},
 \end{equation}
 hence (\ref{lgs}) is satisfied.
Thus we have 
\begin{eqnarray*}
f(x) & = & ({\cal M} +  k\, e^{-\sigma \tau} \sigma) \,  \sin(\omega\, x)  - k \, e^{-\sigma \tau}
\omega
 \,
\cos(\omega\, x),
\\
g(x) & =  &
 - (\sigma + k \, e^{-\sigma \tau}  {\cal M} ) \,  \sin(\omega\, x)   -
 \omega \,
 \cos(\omega\, x)
.
\end{eqnarray*}
The boundary condition (\ref{rb15062018a}) is equivalent to  $f(L) = g(L)$.

{With the choice  $k= - e^{\sigma \tau}$
we have  $f(L) = g(L)$ if $\sigma\in (0,\, {\cal M})$ is such that
}
\[0 = \cos(\sqrt{{\cal M}^2 - \sigma^2}\, L)
.
\]

This is possible if $ {\cal M}\, L$ is sufficiently large in the sense that
\[ {\cal M}  \, L > {\pi} /{2}\]
which follows from assumption (\ref{assumption09082023}).
Note that in this case we have
$|k|  =   e^{\sigma \tau}> 1.$


For $k\not=- e^{\sigma \tau}$ we have  $f(L) = g(L)$ if $\sigma\in (0,\, {\cal M})$ is such that
\begin{equation}
\label{zentral}
H(\sigma, k, \tau):= 
(\sigma + {\cal M}) \, \frac{\tan(\sqrt{{\cal M}^2 - \sigma^2}\, L)}{\sqrt{{\cal M}^2 - \sigma^2}} 
-
\frac { k\, e^{-\sigma \tau}  -1}{ k \, e^{-\sigma \tau} + 1}
=
0
.
\end{equation}

We have 
\[\lim_{\sigma \rightarrow \sqrt{{\cal M}^2 -  \left(\frac{\pi}{ 2\, L} \right)^2} -}
H(\sigma, k, \tau)  =
-  \infty. 
\]

Since
$   \frac{ \sigma + {\cal M} }{  \sqrt{{\cal M}^2 - \sigma^2}}  
=
\sqrt{  \frac{  {\cal M}  +  \sigma  }{ \cal M- \sigma }}    
$
for 
$\sigma_0  \in (0,\, {\cal M})$
and $k\not=- e^{\sigma_0 \tau}$ 
we have 
\[\lim_{\sigma \rightarrow \sigma_0 + }
H(\sigma, k, \tau)  = \sqrt{ \frac{{\cal M} + \sigma_0  }{ {\cal M} - \sigma_0 }} 
\,
{\tan \left( \sqrt{{\cal M}^2 - \sigma^2_0}\, L\right)  }
-  
\frac { k\, e^{-\sigma_0 \tau}  -1}{ k \, e^{-\sigma_0 \tau} + 1}.
\]
Consider the auxiliary function
\[
G_{\sigma_0}(k) = \frac { k\, e^{-\sigma_0 \tau}  -1}{ k \, e^{-\sigma_0 \tau} + 1}.
\]
Then $G_{\sigma_0}$ is continuously differentiable on $( -  e^{\sigma_0 \tau}, \, \infty) $  with the derivative 
\[
G_{\sigma_0}'(k) = \frac{ 2  e^{-\sigma_0 \tau} }{ ( k \, e^{-\sigma_0 \tau} + 1 )^2 }>0.
\]
Thus $G_{\sigma_0}$ is strictly increasing and for
$  k >   -  e^{\sigma_0 \tau}$ with 
$|k| < \hat k$
we have
$G_{\sigma_0}(k) < G_{\sigma_0}( \hat k  )$.
This implies 
\[\lim_{\sigma \rightarrow \sigma_0 + }
H(\sigma, k, \tau)  = \sqrt{ \frac{{\cal M} + \sigma_0  }{ {\cal M} - \sigma_0 }} 
\,
{\tan ( \sqrt{{\cal M}^2 - \sigma^2_0}\, L)  }
-  
G_{\sigma_0}(k)
> 
\sqrt{ \frac{{\cal M} + \sigma_0  }{ {\cal M} - \sigma_0 }} 
\,
{\tan ( \sqrt{{\cal M}^2 - \sigma^2_0}\, L )  }
-  
G_{\sigma_0}(  \hat k  )
=
\]
\[
\sqrt{ \frac{{\cal M} + \sigma_0  }{ {\cal M} - \sigma_0 }} 
\,
{\tan ( \sqrt{{\cal M}^2 - \sigma^2_0}\, L )  }
-  \frac {  \hat k e^{-\sigma_0 \tau}  -1}{  \hat k  e^{-\sigma_0 \tau} + 1}.
\]
Hence if
\begin{equation}
\label{esentialinequality}
\sqrt{ \frac{{\cal M} + \sigma_0  }{ {\cal M} - \sigma_0 }} 
\,
\tan \left( \sqrt{{\cal M}^2 - \sigma^2_0}\, L\right)  
\geq 
\frac{  \hat k  e^{-\sigma_0 \tau}  -1}{  \hat k  e^{-\sigma_0 \tau} + 1}
\end{equation}
the above argument  implies  that  for $  k >   -  e^{\sigma_0 \tau}$ with  $|k| <    \hat k  $ we have
$\lim_{\sigma \rightarrow \sigma_0 + } H(\sigma, k, \tau)  >0$.

Due to (\ref{assumption09082023}) we have
\[  {\tan ( {{\cal M}}  L )}
     \in (-1, \, 0).
     \]
     For $s \in \left[0, \, \sqrt{{\cal M}^2 -  \left(\frac{\pi}{ 2\, L} \right)^2}\right)$
     define  the function
\[
F(s) = \sqrt{ \frac{{\cal M} + s  }{ {\cal M} - s }} 
\, \tan \left( \sqrt{{\cal M}^2 - s^2}\, L\right) . \]
Then $F$  is strictly  decreasing
and 
for $\sigma_0  \in \left (0, \, \sqrt{{\cal M}^2 -  \left(\frac{\pi}{ 2\, L} \right)^2}\right)$
we have
$F(\sigma_0) <   F(0)  = \tan ({\cal M}\, L) <0 $.
Thus the range of $F$ 
for $\sigma_0 \in \left(0, \,   \sqrt{{\cal M}^2 -  \left(\frac{\pi}{ 2\, L} \right)^2}  \right)$
is the interval
$(-\infty,  \,  \tan ({\cal M}\, L) 
) $.
Since $\tan ({\cal M}\, L) + 1 >0$,
we can choose $\sigma_0>0$  sufficiently small such that we have 
$  F(\sigma_0)  +1 >0 $.

Then  we have
\[
\lim_{\tau \rightarrow \infty}
F(\sigma_0) -
\frac{ \hat k  e^{ - \sigma_0 \tau }  -1}{  \hat k e^{ - \sigma_0 \tau  }  + 1} 
=  F(\sigma_0)  +1 >0.
\]

Hence if $ \tau >0$ is sufficiently large we have
\[
F(\sigma_0) -
\frac{  \hat k e^{ - \sigma_0 \tau }  -1}{  \hat k e^{ - \sigma_0 \tau  }  + 1} 
>0.
\]
Thus (\ref{esentialinequality}) holds.
Moreover,  we can choose 
$ \tau >0$  sufficiently large such that 
$   e^{\sigma_0 \tau}< \hat k$.

Therefore for all $\sigma > \sigma_0$ and $|k| <  \hat k  $  we have 
$\lim_{\sigma \rightarrow \sigma_0 + } H(\sigma, k, \tau)  >0$ 
and
\[
\lim_{\sigma \rightarrow \sqrt{{\cal M}^2 -  \left(\frac{\pi}{ 2\, L} \right)^2} -}
H(\sigma, k, \tau)  = -  \infty. 
\]
Hence
due to the continuity of 
$ H(\cdot, k, \tau)$, 
Bolzano's  intermediate value theorem implies  that 
for all $k \in (-\hat k, \, \hat k) $
we can find a number $\sigma \geq \sigma_0$
such that $H(\sigma, k, \tau)  = 0$.
This finishes the proof of   Proposition \ref{proposition1delay}.


%




\section{Hyperbolic weights for quasilinear systems: Stabilization}
\label{quasilinear}

In this section we discuss  how Lyapunov functions
with hyperbolic weights can be used to 
obtain estimates for the stability domains of systems
that are governed by quasilinear PDEs.

We consider a system that is governed by the isothermal Euler equations
with the Riemann invariants $(R_+, R_-)$.
For the stabilization of a stationary state
$(\bar R_+, \, \bar R_-)$
we consider $\delta_\pm = R_\pm - \bar R_\pm$.
We obtain the following system in diagonal form (see \cite{zbMATH05873317}): 
\begin{equation}
    \label{hydrodiag}
\begin{pmatrix}
\delta_+
\\
\delta_-
\end{pmatrix}_t
+
\begin{pmatrix}
\lambda_+(\delta_+, \delta_-) & 0 
\\
 0 & \lambda_-(\delta_+, \delta_-) &
\end{pmatrix}
\begin{pmatrix}
\delta_+
\\
\delta_-
\end{pmatrix}_x
=
\begin{pmatrix}
G_+(\delta_+, \delta_-)
\\
G_-(\delta_+, \delta_-)
\end{pmatrix}
.
\end{equation}
The system is completed with the boundary conditions
\begin{equation}
\label{riemannfeedback0a}
\delta_+(t,\, 0) = k_0 \, \delta_-(t,\, 0)
\end{equation}
and
\begin{equation}
\label{rb15062018astrich}
\delta_-(t,\, L)= k_L \, \delta_+(t,\, L).
\end{equation}
The theory of semi-global solutions (see for example  \cite{li2016exact}, \cite{wang2006exact}) states that 
for a given time horizon $T>0$ 
and $\varepsilon_0 >0$
there exists $\varepsilon_T>0$
such that for all initial states that have a $C^1$-norm
that is less than or equal to $\varepsilon_T>0$
and are $C^1$-compatible with the feedback laws
(\ref{riemannfeedback0a}), (\ref{rb15062018astrich})
 the system has a classical solution on $[0, \, T]$ 
that satisfies the inequalities 
\begin{equation}
\label{sourcetermassumption}
|\partial_x \left( \lambda_\pm( \cdot ) \right)| \leq \varepsilon_0
\end{equation}
and for some $d \geq c>0$ we have
\begin{equation}
\label{aprioribound1lemma}
-d \leq  \lambda_-( \cdot ) \leq -
 c,
 \;\;
 c \leq  \lambda_+( \cdot) \leq d.
\end{equation}
Moreover, for the source term we assume that
\begin{equation}
\label{quelltermvoraussetzung}
|G_\pm(\delta_+, \delta_-)  |\leq {\cal M} \left(|\delta_+| + |\delta_-|\right).
\end{equation}
For  $\psi>0$  and   $\upsilon \in ( \tanh^2(\psi \, L)   , 1]$ 
with the hyperbolic weights that were introduced in Section \ref{defiandprop} 
define the Lyapunov candidate  function
\begin{equation}
    \label{evontdefi}
{\cal E}(t) := 
\frac{1}{2} \int_0^{L}
h_+(x - \tfrac{L}{2} ) \, |\delta_+(t,\, x)|^2 + 
h_-(x - \tfrac{L}{2}  ) \, |\delta_-(t,\, x)|^2 \, dx.
\end{equation}
	For the time-derivative 
of $ {\cal E}(t)  $
 due to (\ref{hydrodiag}) 
 for $t\in [0, \, T]$ we obtain 
	\begin{align*}
		{\cal E}'(t) & =  \int_0^L h_+(x-  \tfrac{L}{2}) \, \delta_+\, (\delta_+)_t + 
   h_-(x-  \tfrac{L}{2})\, \delta_-\, (\delta_-)_t\, dx
		\\
		& = 
		\int_0^L h_+(x-  \tfrac{L}{2}) \, \delta_+\, [- \lambda_+(\delta) \, (\delta_+)_x   +    G_+( \delta ) ] 
		\quad +   h_-(x-  \tfrac{L}{2} ) \, \delta_-\, [-\lambda_-(\delta) \, (\delta_-)_x   +    G_-(\delta) ]\, dx
		\\
		& = 
		\int_0^L    -   h_+(x- \tfrac{L}{2} ) \, \lambda_+(\delta)
                              \,  \left( \tfrac{  \delta_+^2}{2} \right)_x 
		 - 
		h_-(x- \tfrac{L}{2}  )\, \lambda_-(\delta) \,    \left( \tfrac{\delta_-^2}{2} \right)_x 
  \\
	&	+   \left[ \delta_+ \,G_+(\delta) \, h_+(x-  \tfrac{L}{2} ) + \delta_- \,  G_-(\delta) \, h_-(x-  \tfrac{L}{2} )\right]   \, dx.
	\end{align*}
 Now (\ref{quelltermvoraussetzung})
 and
 $|\delta_+\, \delta_-|  \leq \frac{|\delta_+|^2}{2} + \frac{|\delta_-|^2}{2}$
 yield  the inequality
\[
 \int_0^L \left[ \delta_+ \,G_+(\delta) \, h_+(x-  \tfrac{L}{2} ) + \delta_- \,  G_-(\delta) \, h_-(x-  \tfrac{L}{2} )\right]   \, dx
 \]
 \[
 \leq \int_0^L {\cal M} \,  h_+(x-  \tfrac{L}{2})  \left( \tfrac{3}{2} \delta_+^2 + \tfrac{1}{2} \delta_-^2  \right) 
   + {\cal M} \,  h_-(x-  \tfrac{L}{2})  \left( \tfrac{1}{2} \delta_+^2 + \tfrac{3}{2} \delta_-^2  \right) \, dx.
 \]
With integration by parts and Lemma \ref{lemma1} e) the above inequalities  yield
	\begin{align*}
	        {\cal E}'(t)     = 
		&  \int_0^L  
   \left[  h_+'(x-  \tfrac{L}{2} ) \, \lambda_+(\delta)  +  h_+(x-  \tfrac{L}{2}  )  \, \partial_x \lambda_+(\delta) \right] \tfrac{\delta_+^2 }{2} 
    \\
    & +    \left[  h_-'(x-   \tfrac{L}{2}) \, \lambda_-(\delta)  +  h_-(x-  \tfrac{L}{2} )  \, \partial_x \lambda_-(\delta) \right] \tfrac{\delta_-^2 }{2} 
  	 dx
     \\
		&+
		\left[ -   h_+(x- \tfrac{L}{2} ) \, \lambda_+(\delta)
                              \,   \tfrac{  \delta_+^2}{2} 
		 - 
		h_-(x- \tfrac{L}{2}  )\, \lambda_-(\delta) \,     \tfrac{\delta_-^2}{2} 
  \Big|_{x=0}^{x=L} \right]
  \\
	&	+   \int_0^L \left[ \delta_+ \,G_+(\delta) \, h_+(x-  \tfrac{L}{2} ) + \delta_- \,  G_-(\delta) \, h_-(x-  \tfrac{L}{2} )\right]   \, dx
  \\
  &  \leq - \tfrac{\psi}{2} \int_0^L \left[  \frac{  1 +  \upsilon   }{  2 \sqrt{\upsilon}   }      \right] h_+(x-  \tfrac{L}{2}  )\,\lambda_+(\delta) \, \delta_+^2
  + \left[ \frac{ 1 - \upsilon   }{   2 \sqrt{\upsilon}     }  \right]   h_-(x-  \tfrac{L}{2} ) \lambda_+(\delta) \,\delta_+^2
   \\ &  + \left[ \frac{ 1 - \upsilon  }{   2 \sqrt{\upsilon}     }  \right] 
 h_+(x  -  \tfrac{L}{2} ) \, | \lambda_-(\delta)|  \, \delta_-^2 
  + \left[  \frac{   1 + \upsilon  }{  2 \sqrt{\upsilon}   }      \right]
 h_-(x-   \tfrac{L}{2} ) \,  |\lambda_-(\delta)| \, \delta_-^2 \, dx 
   \\
	&	+   \int_0^L   h_+(x-  \tfrac{L}{2}  )  \, \partial_x \lambda_+(\delta) \,  \tfrac{\delta_+^2 }{2}  +    
                           h_-(x-  \tfrac{L}{2}  )  \, \partial_x \lambda_-(\delta) \,  \tfrac{\delta_-^2 }{2}    \, dx
   \\
	& +	\int_0^L {\cal M} \,  h_+(x-  \tfrac{L}{2})  \left( \tfrac{3}{2} \delta_+^2 + \tfrac{1}{2} \delta_-^2  \right) 
   + {\cal M} \,  h_-(x-  \tfrac{L}{2})  \left( \tfrac{1}{2} \delta_+^2 + \tfrac{3}{2} \delta_-^2  \right) \, dx
 \\
 	&+
		\left[ -   h_+(x- \tfrac{L}{2} ) \, \lambda_+(\delta)
                              \,   \tfrac{  \delta_+^2}{2} 
		 - 
		h_-(x- \tfrac{L}{2}  )\, \lambda_-(\delta) \,     \tfrac{\delta_-^2}{2} 
  \Big|_{x=0}^{x=L} \right].
	\end{align*}

 If the boundary term in the last line is less than or equal to zero, 
  since  $\upsilon \in ( \tanh^2(\psi \, L)   , 1]$ due to Lemma  \ref {lemma1}
  and (\ref{sourcetermassumption})   we obtain the inequality 
\begin{align*}
    {\cal E}'(t)    \leq 
  & - \tfrac{\psi}{2}  \int_0^L \left[ 
  c\,
  \frac{  1 + \upsilon   }{  2 \sqrt{\upsilon}   }  -   3 \frac{{\cal M}}{\psi}  -    \frac{{\varepsilon_0 }}{\psi}
  \right] h_+(x -  \tfrac{L}{2}  )\, \delta_+^2
  + \left[  
  c\, \frac{ 1 - \upsilon   }{   2 \sqrt{\upsilon}     } 
   -        \frac{{\cal M}}{\psi}
  \right] h_-(x  -   \tfrac{L}{2} ) \,  \delta_+^2
   \\ &  + \left[  
   c\, \frac{ 1 - \upsilon  }{   2 \sqrt{\upsilon}     }  -    \frac{{\cal M}}{\psi}  \right] 
 h_+(x-  \tfrac{L}{2})  \, (\delta_-)^2 
  + \left[ 
  c\, \frac{ 1 +   \upsilon   }{  2 \sqrt{\upsilon}   }        -   3 \frac{{\cal M}}{\psi}  -    \frac{{\varepsilon_0 }}{\psi} \right]
 h_-(x-  \tfrac{L}{2})   \, (\delta_-)^2 \, dx. 
	\end{align*}
 
If
\begin{equation}
\label{1007202312asb}
3 \, {\cal M}  + \varepsilon_0 
\leq
 c\, 
\psi \, 
\frac{ 1 - \upsilon  }{  2\,  \sqrt{\upsilon}     }
=
\frac{c}{2} 
\, {\psi} 
\left(
\frac{1}{    \sqrt{\upsilon}     }
-
{    \sqrt{\upsilon}     }
\right)
\end{equation}
we have
\[
\frac{{\cal M}}{\psi} 
\leq
 c\, 
\frac{ 1 - \upsilon  }{  2\,  \sqrt{\upsilon}     }
\]
and  this yields
\begin{align*}
    {\cal E}'(t)    \leq  - 
    \psi \,
    \frac{c}{2} 
\left[  \frac{  1 + \upsilon   }{  2 \sqrt{\upsilon}   }   -    \frac{ 1 - \upsilon  }{   2 \sqrt{\upsilon}     } \right]   
\int_0^{L}
h_+(x - \tfrac{L}{2} ) \, |\delta_+(t,\, x)|^2 + 
h_-(x - \tfrac{L}{2}  ) \, |\delta_-(t,\, x)|^2 \, dx
=
 -
c\,
 \psi 
\,
\sqrt{\upsilon} 
\,
{\cal E}(t).
	\end{align*}
	By Gronwall's inequality this implies
	$ {\cal E}(t) \leq \exp(-   
 c 
 \,\psi 
\,
\sqrt{\upsilon}  \, t) \, 
{\cal E}(0),
$
	so it remains to check the assumption on
 the boundary term.
 Below we will derive sufficient conditions that are stated in  (\ref{suffsmallgain}).

 Note that (\ref{1007202312asb} ) is equivalent to
\begin{equation}
\label{1007202312}
3 {\cal M} L + \varepsilon_0 L
\leq 
  \frac{c}{2 \,} 
\psi \, L
\frac{ 1 - \upsilon  }{    \sqrt{\upsilon}     }
= 
  \frac{c}{2} \,
\psi L
\left(
\frac{1}{    \sqrt{\upsilon}    }
-
{    \sqrt{\upsilon}     }
\right).
\end{equation}
For the right-hand-side of
(\ref{1007202312})  we have  the  upper  bound
\[
 \frac{c}{2} 
 \sup_{\psi >0 }
\sup_{ \upsilon \in (\tanh^2(\psi  \, L), 1]   }
 \psi L
\left(
\frac{1}{    \sqrt{\upsilon}    }
-
{    \sqrt{\upsilon}     }
\right)
=
 \frac{c}{2} 
\sup_{ z >0}
 z
\left(
\frac{1}{    \tanh(z) }
-
{   \tanh(z)    }
\right)
=
 \frac{c}{2} 
.
\]
(Note that the value
$ \frac{c}{2}  $ is attained in the limit $z \rightarrow 0$.)

So  we obtain an upper bound for
the values of  ${\cal M} L$ for which we can guarantee that
the system is exponentially stable, namely

\[
3 {\cal M} L + \varepsilon_0 L 
 < 
 \frac{c}{2} 
.
\]
This bound is  useful in the application in gas dynamics
where $c$ and $d$
are  related  to the sound speed that is quite large, see 
\cite{zbMATH05142458}.

The condition on the boundary terms
\[  h_+(- \tfrac{L}{2} ) \, \lambda_+(\delta( 0  )  )
                              \,   \tfrac{  k_0^2 \delta_-^2( 0) }{2} 
		 +
		h_-(- \tfrac{L}{2}  )\, \lambda_-(\delta ( 0 )  ) \,     \tfrac{\delta_-^2( 0 ) }{2} 
 -   h_+( \tfrac{L}{2} ) \, \lambda_+(\delta( L ) )
                              \,   \tfrac{  \delta_+^2( L )  }{2} 
		 - 
		h_-( \tfrac{L}{2}  )\, \lambda_-(\delta( L ) ) \,     \tfrac{  k_L^2  \delta_+^2(  L )}{2} 
    \leq 0
  \]
is satisfied if
\[
h_+(- \tfrac{L}{2} ) \, \lambda_+(\delta( 0  )  )
                              \,    k_0^2  \leq - 	h_-(- \tfrac{L}{2}  )\, \lambda_-(\delta ( 0 )  ) 
\]
and
\[
- 
		h_-( \tfrac{L}{2}  )\, \lambda_-(\delta( L ) ) \,  k_L^2
  \leq
   h_+( \tfrac{L}{2} ) \, \lambda_+(\delta( L ) ).
\]
This is equivalent to
  \[
 k_0^2  \leq  \frac{  h_-(- \tfrac{L}{2}  )}{ h_+(- \tfrac{L}{2} )  } \, \frac{  |\lambda_-(\delta ( 0 )  ) |}{\lambda_+(\delta( 0  )  )}
 =
\frac{  h_+(\tfrac{L}{2}  )}{ h_-( \tfrac{L}{2} )  } \, \frac{  |\lambda_-(\delta ( 0 )  ) |}{\lambda_+(\delta( 0  )  )}
 ,
  \;\;
   k_L^2  \leq  \frac{  h_+( \tfrac{L}{2}  )}{ h_-( \tfrac{L}{2} )  } \, \frac{  \lambda_+(\delta ( L )  ) }{ |\lambda_-(\delta( L  ) | )}
   .
\]
Sufficient conditions are
(with $i\in \{0,\, L\}$)
  \begin{equation}
  \label{suffsmallgain}
 k_i^2  \leq 
 \, 
\frac{c}{d}
\frac{  h_+(\tfrac{L}{2}  )}{ h_-( \tfrac{L}{2} )  } 
=
\frac{c}{d}
\frac{ \sqrt{\upsilon} \left(1 + \tfrac{1}{\cosh(\psi\, L)} \right) - \tanh(\psi \, L)}{ \sqrt{\upsilon} \left(1 + \tfrac{1}{\cosh(\psi\, L)} \right) +  \tanh(\psi \, L)}.
\end{equation}
where the last equality follows from  (\ref{quationt23022024}).


Thus we have shown the following results for the boundary control of
the quasilinear system:
\begin{theorem}
	\label{proposition2223}
Assume that the source term $G_\pm$  satisfies (\ref{quelltermvoraussetzung}). 
 
 Assume that  the	closed loop system
	(\ref{hydrodiag}),  
	(\ref{riemannfeedback0a}),
	(\ref{rb15062018astrich})
 has a classical solution on $[0, \, T]$ such that 
    (\ref{sourcetermassumption}) 
    and
  (\ref{aprioribound1lemma}) 
hold for all $t \in [0,\, T]$, that we have 
\begin{equation}
    \label{lbedingung}
 \left( 2 \, \varepsilon_0 + 6\, {\cal M} \right) \, L 
 < 
c 
\end{equation}
        and the feedback gains 
	$|k_0|$ and $|k_L|$ are  sufficiently small
 (for example, that  (\ref{suffsmallgain}) holds). 

 Then 	the	closed loop system
	(\ref{hydrodiag}),  
	(\ref{riemannfeedback0a}),
	(\ref{rb15062018astrich})
 decays exponentially fast on
 $[0, \, T]$.
If   $\psi>0$  and   $\upsilon \in ( \tanh^2(\psi \, L)   , 1]$ 
are chosen such that (\ref{1007202312}) holds, 
 the Lyapunov function
$ {\cal E}(t)$  decays exponentially  with the rate
\[     c \,\psi \,\sqrt{\upsilon} . \]

\end{theorem}

\textbf{Remark:}
In order to make sure that the system is globally well-posed
for the time intervall $[0, \,\infty)$ 
also Lyapunov functions for
the first and second derivatives can be considered.
This yields the exponential decay of solutions
with values in $H^2(0, L)$,
see \cite{zbMATH06742239}, \cite{hayat2021exponential}. 
Solutions in $H^2$ are also studied in \cite{MR3561145}.

In order to show the exponential decay of an $H^1$-Lyapunov function,
we have to assume that $G$ has partial derivatives that 
are continuous and 
an additional assumption for the derivatives of the source term is necessary:
\begin{equation}
\label{ableitungsvor}
| \partial_\pm G_\pm(\delta_+, \delta_-)  |\leq {\cal M} \left(|\delta_+| + |\delta_-|\right).
\end{equation}

Similarly, in order to show the exponential decay of an $H^2$-Lyapunov function,
we have to assume that $G$ has second order  partial derivatives that 
are continuous and 
an additional assumption for the second derivatives of the source term is necessary:
\begin{equation}
\label{ableitungsvor2}
z^\top
\begin{pmatrix}
\partial_{++} G_\pm & \partial_{+-} G_\pm 
\\
\partial_{+-} G_\pm & \partial_{--} G_\pm 
\end{pmatrix}
z \leq {\cal M}
.
\end{equation}



\textbf{Remark:}
Condition  (\ref{lbedingung}) requires that $L$ 
is sufficiently small.
Since often in the applications the right hand side  $c$  
can be chosen  proportional to the sound speed,
 (\ref{lbedingung}) is valid for interesting lengths.
The bound  $\varepsilon_0$ can often be chosen quite small in
the applications since only small changes in the states occur.
Horizontal pipes with anti-fricion coating allow for small values of ${\cal M}$.

Compared with Lemma 5.2. in \cite{zbMATH05873317},
conditions (\ref{lbedingung}), (\ref{1007202312asb}) and  (\ref{suffsmallgain})  have the advantage that they 
can be verified more  easily.




\section{Hyperbolic weights for quasilinear systems: Instability}
\label{quasilinearin}
In this section we show that for sufficiently large
values of $L$, boundary feedback stabilization 
of quasilinear hyperbolic systems in general is not possible.

We assume that there exists a number $  {\cal N}>\varepsilon_0>0$  such that
for the source term, we have 
\begin{equation}
\label{quelltermvoraussetzungdesta}
\delta_\pm \,  G_\pm(\delta_+, \delta_-)  \geq   
{\cal N} |\delta_\pm|^2.
\end{equation}

Then we have
\begin{eqnarray*}
    {\cal E}'(t)  
    &
\geq
&
 - \tfrac{\psi}{2} \int_0^L \left[  \frac{  1 +  \upsilon   }{  2 \sqrt{\upsilon}   }      \right] h_+(x-  \tfrac{L}{2}  )\,\lambda_+(\delta) \, \delta_+^2
  + \left[ \frac{ 1 - \upsilon   }{   2 \sqrt{\upsilon}     }  \right]   h_-(x-  \tfrac{L}{2} ) \lambda_+(\delta) \,\delta_+^2
   \\
   &  + &
   \left[ \frac{ 1 - \upsilon  }{   2 \sqrt{\upsilon}     }  \right] 
 h_+(x  -  \tfrac{L}{2} ) \,  |\lambda_-(\delta)|  \, \delta_-^2 
  + \left[  \frac{   1 + \upsilon  }{  2 \sqrt{\upsilon}   }      \right]
 h_-(x-   \tfrac{L}{2} ) \,  |\lambda_-(\delta)| \, \delta_-^2 \, dx 
   \\
	&	+  &  \int_0^L   h_+(x-  \tfrac{L}{2}  )  \, \partial_x \lambda_+(\delta) \,  \tfrac{\delta_+^2 }{2}  +    
                           h_-(x-  \tfrac{L}{2}  )  \, \partial_x \lambda_-(\delta) \,  \tfrac{\delta_-^2 }{2}    \, dx
   \\
	& + &	\int_0^L {\cal N} \,  h_+(x-  \tfrac{L}{2})  \,  \delta_+^2 
   + {\cal N} \,  h_-(x-  \tfrac{L}{2})  \, \delta_-^2  \, dx
 \\
 	& + & 
		\left[ -   h_+(x- \tfrac{L}{2} ) \, \lambda_+(\delta)
                              \,   \tfrac{  \delta_+^2}{2} 
		 - 
		h_-(x- \tfrac{L}{2}  )\, \lambda_-(\delta) \,     \tfrac{\delta_-^2}{2} 
  \Big|_{x=0}^{x=L} \right].
	\end{eqnarray*}
{\bl  With the choice $\upsilon = 1$}  
 due to Lemma  \ref {lemma1}
  and (\ref{sourcetermassumption})   we obtain the inequality 
\begin{align*}
    {\cal E}'(t)    \geq 
  & - \tfrac{\psi}{2}  \int_0^L \left[ 
  c\,
  \frac{  1 + \upsilon   }{  2 \sqrt{\upsilon}   }  -    \frac{{\cal N}}{\psi}  +  \frac{{\varepsilon_0 }}{\psi}
  \right] h_+(x -  \tfrac{L}{2}  )\, \delta_+^2
  + \left[  
  c\, \frac{ 1 - \upsilon   }{   2 \sqrt{\upsilon}     } 
  \right] h_-(x  -   \tfrac{L}{2} ) \,  \delta_+^2
   \\ &  + \left[  
   c\, \frac{ 1 - \upsilon  }{   2 \sqrt{\upsilon}     }   \right] 
 h_+(x-  \tfrac{L}{2})  \, (\delta_-)^2 
  + \left[ 
  c\, \frac{ 1 +   \upsilon   }{  2 \sqrt{\upsilon}   }        -    \frac{{\cal N}}{\psi} +   \frac{{\varepsilon_0 }}{\psi} \right]
 h_-(x-  \tfrac{L}{2})   \, (\delta_-)^2 \, dx
 \\
 	&+
		\left[ -   h_+(x- \tfrac{L}{2} ) \, \lambda_+(\delta)
                              \,   \tfrac{  \delta_+^2}{2} 
		 - 
		h_-(x- \tfrac{L}{2}  )\, \lambda_-(\delta) \,     \tfrac{\delta_-^2}{2} 
  \Big|_{x=0}^{x=L} \right]
  \\
   \geq &
    \tfrac{\psi}{2}  \int_0^L \left[ 
 \frac{{\cal N}}{\psi}  -   \frac{{\varepsilon_0 }}{\psi}
    -
  c\,
  \frac{  1 + \upsilon   }{  2 \sqrt{\upsilon}   }  
  \right] h_+(x -  \tfrac{L}{2}  )\, \delta_+^2
      +
  \left[
   \frac{{\cal N}}{\psi}  -   \frac{{\varepsilon_0 }}{\psi} -
  c\, \frac{ 1 +   \upsilon   }{  2 \sqrt{\upsilon}   }       \right]
 h_-(x-  \tfrac{L}{2})   \, (\delta_-)^2 \, dx
 \end{align*}
\begin{equation}
\label{26022024a}
 + 
		\left[ -   h_+(x- \tfrac{L}{2} ) \, \lambda_+(\delta)
                              \,   \tfrac{  \delta_+^2}{2} 
		 - 
		h_-(x- \tfrac{L}{2}  )\, \lambda_-(\delta) \,     \tfrac{\delta_-^2}{2} 
  \Big|_{x=0}^{x=L} \right]
  .
  \end{equation}
For 
$\psi>0$  
{\bl we  study  the inequality} 
\begin{equation}
\label{2602a}
- \tfrac{\psi}{2}  \left[ c\,  \frac{  1 +  \upsilon   }{  2 \sqrt{\upsilon}   }   
 -    \frac{{\cal N}}{\psi}  +   \frac{{\varepsilon_0 }}{\psi} \right]
>0.
\end{equation}


Let $\eta>0$ be given.
 With the choice 
$\psi = \frac{\eta}{L}$
due to $\upsilon = 1$ 
(\ref{2602a}) 
yields
\begin{equation}
\label{2602b1}
 \tfrac{\eta}{2\, L}  \left[
   \frac{L}{\eta}  ({\cal N} -   {{\varepsilon_0 }}) -
 c
 \right]
>0.
\end{equation}
Inequality (\ref{2602b1}) is equivalent to
\begin{equation}
\label{2602c1}
 ({\cal N} -   {{\varepsilon_0 }})  L
>
\eta\,  c.
\end{equation}




The last  term in (\ref{26022024a})
is greater  than or equal to zero
if
\[
h_+(- \tfrac{L}{2} ) \, \lambda_+(\delta( 0  )  )
                              \,    k_0^2  \geq - 	h_-(- \tfrac{L}{2}  )\, \lambda_-(\delta ( 0 )  ) 
\]
and
\[
- 
		h_-( \tfrac{L}{2}  )\, \lambda_-(\delta( L ) ) \,  k_L^2
  \geq
   h_+( \tfrac{L}{2} ) \, \lambda_+(\delta( L ) ).
\]
This is equivalent to
  \[
 k_0^2  \geq  \frac{  h_-(- \tfrac{L}{2}  )}{ h_+(- \tfrac{L}{2} )  } \, \frac{  |\lambda_-(\delta ( 0 )  ) |}{\lambda_+(\delta( 0  )  )}
 =
\frac{  h_+(\tfrac{L}{2}  )}{ h_-( \tfrac{L}{2} )  } \, \frac{  |\lambda_-(\delta ( 0 )  ) |}{\lambda_+(\delta( 0  )  )}
 ,
  \;\;
   k_L^2  \geq  \frac{  h_+( \tfrac{L}{2}  )}{ h_-( \tfrac{L}{2} )  } \, \frac{  \lambda_+(\delta ( L )  ) }{ |\lambda_-(\delta( L  ) | )}
   .
\]
Sufficient conditions are
(with $i\in \{0,\, L\}$)
  \begin{equation}
  \label{suffsmallgainin}
 k_i^2  \geq 
 \, 
\frac{d}{c}
\frac{  h_+(\tfrac{L}{2}  )}{ h_-( \tfrac{L}{2} )  } 
=
\frac{d}{c}
\frac{ \sqrt{\upsilon} \left(1 + \tfrac{1}{\cosh(\psi\, L)} \right) - \tanh(\psi \, L)}{ \sqrt{\upsilon} \left(1 + \tfrac{1}{\cosh(\psi\, L)} \right) +  \tanh(\psi \, L)}.
\end{equation}
where the last equation follows from  (\ref{quationt23022024}).
Note that we have
\[
\lim_{L\rightarrow\infty, \,  \upsilon \in ( \tanh^2(\psi \, L)   , 1]}
\frac{ \sqrt{\upsilon} \left(1 + \tfrac{1}{\cosh(\psi\, L)} \right) - \tanh(\psi \, L)}{ \sqrt{\upsilon} \left(1 + \tfrac{1}{\cosh(\psi\, L)} \right) +  \tanh(\psi \, L)}
=0.
\]

In particular for $\psi = \frac{\eta}{L}$ and $ \upsilon  = 1   $
inequality (\ref{suffsmallgainin}) yields  the sufficient condition
\[k_i^2  \geq \frac{d}{c} \exp\left(- \eta \right).\]

Hence for all 
$k_0 \not=0 $ 
and  
$k_L \not =0$ 
if $L$ is sufficiently large  we obtain the inequality 
\begin{align*}
    {\cal E}'(t)    \geq 0.
	\end{align*}
Thus  ${\cal E}(t) $   
does not decrease 
and thus the system is not asymptotically stable.
Thus we have shown the following Proposition about the instability
of the quasilinear system for large values of $L$:

\begin{theorem}
	\label{proposition2223in}

 Assume that  the	closed loop system
	(\ref{hydrodiag}),  
	(\ref{riemannfeedback0a}),
	(\ref{rb15062018astrich})
 has a classical solution on $[0, \, T]$ such that 
    (\ref{sourcetermassumption}) 
    and
  (\ref{aprioribound1lemma}) 
hold for all $t \in [0,\, T]$
        and the feedback gains  satisfy
	$|k_0|\not=0$ and $|k_L|\not=0$.
  We assume that there exists a number $  {\cal N}>\varepsilon_0>0$  such that
the  source term satisfies  (\ref{quelltermvoraussetzungdesta}).

 Then if $L$ is sufficiently large 	the	closed loop system
	(\ref{hydrodiag}),  
	(\ref{riemannfeedback0a}),
	(\ref{rb15062018astrich})
 does not decay  on  $[0, \, T]$.

 To be precise,  if   $\psi>0$  and   $\upsilon = 1$ 
are chosen such that 
(\ref{2602a}) and 
(\ref{suffsmallgainin}) hold, 
 the Lyapunov function
$ {\cal E}(t)$  
(see \ref{evontdefi}) 
does not decrease 
on $[0, T]$.
%
%
%
 A sufficient condition for instability 
   with a real parameter $\eta>0$  is
\begin{equation}
\label{2602cordnung}
 ({\cal N} -   {{\varepsilon_0 }})  \,  L
>
\eta\, c,
 \;
 |k_\iota|^2  \geq 
\frac{d}{c}
\,
 \exp( - \eta),\;  (\iota \in \{0,\, L\})
.
\end{equation}




\end{theorem}
\begin{remark}
In the  sufficient condition (\ref{2602cordnung}) 
the left-hand side 
of the first inequality 
grows linearly with $L$,
similarly as in (\ref{lbedingung}).


\end{remark}


\section{Conclusions}
We have introduced hyperbolic weights for
quadratic Lypunov functions. We have shown that 
in certain cases they yield larger  stability
regions than exponential weights.
We have studied this in the context of stabilization problems.
We expect that the hyperbolic weights
are also useful in the analysis of the synchronization of
observers, see \cite{gugat2021exponential}.
It would be very interesting to investigate whether the
hyperbolic weights are also suitable for
the analysis of networked systems,
see for example \cite{zbMATH07117070}.

{\bf Acknowledgments.}
We are grateful for  support by DFG
in the framework of the Collaborative Research Centre
CRC/Transregio 154,
Mathematical Modelling, Simulation and Optimization Using the Example of Gas Networks,
project C03 and  project  C05, Projektnummer 239904186
%
and
the Bundesministerium für Bildung und Forschung (BMBF)
and the Croatian Ministry of Science and Education 
under DAAD grant 57654073 'Uncertain data in control of PDE systems'.


\bibliographystyle{alpha}
\bibliography{main}

\end{document}